\documentclass[12pt]{article}
\usepackage{amssymb}
\usepackage{amsmath}

\usepackage{latexsym}
\usepackage{amsfonts}
\usepackage{times}

\newtheorem{theorem}{Theorem}[section]
\newtheorem{lemma}[theorem]{Lemma}

\newenvironment{Proof}{\removelastskip\par\medskip
\noindent{\em Proof.} \rm}{\penalty-20\null\hfill$\square$\par\medbreak}

\begin{document}

\centerline{\textbf{\Large Inequalities for convolutions of functions }}

\centerline{\textbf{\Large  on commutative hypergroups}}

\vspace{3mm} \centerline{by}

\vspace{3mm} \centerline{\textbf{Mubariz G. Hajibayov}}

 \centerline{\it National Aviation Academy, Baku, Azerbaijan}
\centerline{\it and}
 \centerline{\it Institute of Mathematics and Mechanics, Baku, Azerbaijan}

\centerline{\it (hajibayovm@yahoo.com)}

\begin{abstract}

The generalized Young inequality on the Lorentz spaces for commutative hypergroups  is introdused and  an application of it is given to the theory of fractional integrals. The
boundedness on the Lorentz space and  the Hardy-Littlewood-Sobolev theorem for the fractional integrals on the commutative hypergroups is proved.
\end{abstract}

\vskip 5pt


\vskip 5pt


\section{Introduction and preliminaries}

It is known that a convolution of two functions on $R^n$ is defined by
 $$
 f\ast_{R^n} g(x)=\int \limits_{R^n}f(x-y)g(y)dy.
 $$

Classical Young's inequality on the $L^p(R^n)$ spaces the convolution of two functions on $R^n$  states that
if $f \in L^p(R^n)$ and $g \in L^q(R^n)$, then
$$
\|f\ast_{R^n} g\|_r\leq C\|f\|_p\| g\|_q,
$$
 where $p,q \in [1,\infty]$ and $\frac 1p + \frac 1q=\frac 1r+1$.

 The generalized Young inequality give us the boundedness on the Lorentz spaces for the convolution of two functions on $R^n$.
 \begin{theorem}( \cite{Z} Theorem 2.10.1 )
  If  $f\in L^{p_1,q_1} \left( R^n \right)$, $\varphi \in L^{p_2,q_2} \left( R^n \right)$ and $\dfrac{1}{p_1}+\dfrac{1}{p_2}>1$, then $(f\ast\varphi)\in L^{p_0,q_0} \left( R^n \right)$ where $\dfrac{1}{p_1}+\dfrac{1}{p_2}-1=\dfrac{1}{p_0}$ and $q_0\geq 1$ is any number such that
  $\dfrac{1}{q_1}+\dfrac{1}{q_2}\geq \dfrac{1}{q_0}$.\\
  Moreover,
  $$
\|(f\ast\varphi)\|_{K,p_0,q_0}\leq 3p_0\|f\|_{K,p_1,q_1}\|\varphi\|_{K,p_2,q_2}.
$$
    \end{theorem}
An extension of the Young inequality to the convolution
$$
f\ast_G g(x)=\int _{G}f(xy)g(y^{-1})d\mu(y),
$$
where $\mu$ is the Haar measure on local compact group $G$, was given in \cite{HR}(see Theorem 20.18 in \cite{HR}  ).


In the theory of locally compact groups there arise certain spaces which, though not groups, have some of the structure of groups. Often, the structure can be expressed in terms of an abstract convolution of measures on the space. 

A hypergroup  $(K,\ast_K)$ consists of a locally compact Hausdorff space K together with a
bilinear, associative, weakly continuous convolution on the Banach space of all bounded
regular Borel measures on $K$ with the following properties:
\begin{itemize}
\item[1.] For all $x,y \in K$, the convolution of the point measures $\delta _x\ast_K \delta _y$ is a probability measure with compact support.

\item[2.] The mapping: $K\times K\rightarrow \mathcal{C}(K)$, $(x,y)\mapsto\, \text{supp}\,(\delta _x\ast_K \delta _y)$ is continuous with respect to the Michael topology on the space $\mathcal{C}(K)$ of all nonvoid compact subsets of $K$, where this topology
is generated by the sets
$$U_{V,W}=\{L\in\mathcal{C}(K): L\cap V\neq \varnothing, L\subset W \}$$
with $V,W$ open in $K$.

\item[3.] There is an identity $e \in K$ with
$\delta _e\ast_K \delta _x=\delta _x\ast_K \delta _e=\delta _x$ for all $x \in K$.

\item[4.] There is a continuous involution $\thicksim$ on $K$ such that
$$
\left(\delta _x\ast_K \delta _y\right)^\thicksim =\delta _{y^\thicksim }\ast_K \delta _{x^\thicksim}
$$
and
$e \in supp(\delta _x\ast_K \delta _y) \Leftrightarrow x=y^\thicksim$ for $x,y \in K$
(see \cite{J}, \cite{Sp}, \cite{BH}).
\end{itemize}
A hypergroup $K$ is called commutative if $\delta _x\ast_K \delta _y=\delta _y\ast_K \delta _x$  for all
$x, y\in K$. It is well known that every commutative hypergroup $K$ possesses a Haar measure which will be denoted by $\lambda$
(see \cite{Sp}). That is, for every Borel measurable function $f$ on $K$,
$$
\int \limits_K f(\delta _x\ast_K \delta _y)d\lambda (y)=\int \limits_K f(y)d\lambda (y) \,\,\, (x \in K).
$$
Define the generalized translation operators $T^x$, $x\in K$,
by
$$
T^xf(y) =\int \limits_Kfd(\delta _x\ast_K \delta _y)
$$
for all $y \in K$. If $K$ is a commutative hypergroup, then $T^xf(y)=T^yf(x)$ and
the convolution of two functions is defined by
$$
f\ast_K \varphi(x)=\int \limits_KT^xf(y)\varphi(y^\thicksim)d\lambda (y).
$$
Note that $f\ast_K \varphi=\varphi\ast_K f $.\\
For $1\leq p\leq \infty$, the Lebesgue space $L^{p} \left( K,\lambda \right) $ is defined as
 $$
 L^{p} \left( K,\lambda \right)=\{f: f \,\text{is} \,\lambda \, \text{-measurable on}\, K, \|f\|_{K,p}<\infty \}
 $$
 where $\|f\|_{K,p}$ is defined by
 $$
 \|f\|_{K,p}=\begin{cases}
 \left(
\int\limits_{K}\left| f\left( x\right) \right| ^{p} d\lambda \left(
x\right) \right) ^{\frac{1}{p} },&\text{if}\,\, 1\leq p<\infty\\
 \text{ess}\sup \limits _{x\in K}f(x), &  \text{if}\,\, p= \infty.
 \end{cases}
 $$
Let $1\leq p\leq \infty$. If $f$ is in $L^{p} \left( K,\lambda \right) $ and $\varphi$ is in $L^{1} \left( K,\lambda \right) $, then the function $f\ast_K \varphi$
belongs to $L^{p} \left( K,\lambda \right) $ and
$$
\|f\ast_K \varphi\|_{K,p}\leq \|f\|_{K,p}\| \varphi\|_{K,1}
$$
Let $f$ be a $\lambda$-measurable function defined on the hypergroup $K$. The distribution function
$\lambda _f$ of the function $f$  is given by
$$
\lambda _f(s)=\lambda \{x:x \in K, |f(x)|>s\}, \,\,\text{for} \,\,s\geq 0.
$$
The distribution function $\lambda _f$ is non-negative, non-increasing and continuous from the right. With the distribution function we associate the non-increasing rearrangement of $f$ on $[0,\infty )$ defined by
$$
f^{\ast_K}(t)=\inf \{s>0:\lambda _f(s) \leq t\}.
$$
 Some elementary properties of $\lambda _f$ and $f^{\ast_K}$  are listed below. The proofs of them can be found in \cite{BSh}.
 \begin{itemize}
 \item[$(1)$]
  If $\lambda _f$ is continuous and strictly decreasing, then $f^{\ast_K}$ is the inverse of $\lambda _f$, that is $f^{\ast_K}=\left(\lambda _f\right)^{-1}$.
 \item[$(2)$]
 $f^{\ast_K}$ is continuous from the right.
  \item[$(3)$]
 $$
 m _{f^{\ast_K}}(s)=\lambda _f(s),\,\,\text{ for  all}\, \, s>0,
 $$
 where $m _{f^{\ast_K}}$ is a distribution function
of the function $f^{\ast_K}$ with respect to Lebesgue measure $m$ on $(0,\infty)$.
  \item[$(4)$]
  \begin{equation} \label{fstar}
  \int\limits_0^t f^{\ast_K}(s)ds=tf^{\ast_K}(t)+\int\limits_{f^{\ast_K}(t)}^\infty \lambda _f(s)ds
  \end{equation}
  \item[$(5)$]
 If $f \in L^{p} \left( K,\lambda \right)$, $1\leq p< \infty$, then
 $$
 \left(\int\limits_{K}\left| f\left( x\right) \right| ^{p} d\lambda \left(
x\right) \right) ^{\frac{1}{p} }=\left(p\int \limits _0^\infty s^{p-1}\lambda _f(s)ds\right)^{\frac 1p}=\left(\int \limits _0^\infty \left(f^{\ast_K}(t)\right)^pdt\right)^{\frac 1p}.
 $$
 Furthermore, in the case $p=\infty$,
 $$
 \text{ess}\sup \limits_{x\in K}f(x)=\inf \{s: \lambda_f (s)=0\}=f^{\ast_K} (0)
 $$
 \end{itemize}
  $f^{\ast \ast_K}$ will denote the maximal function of $f^{\ast_K }$ defined by
  $$
 f^{\ast \ast_K}(t)=\frac 1t \int \limits _0^t f^{\ast_K}(u)du, \,\,\text{for}\,\,t>0.
 $$
 Note the following properties of $f^{\ast \ast_K}$:
 \begin{itemize}
 \item[$(1')$]
 $f^{\ast \ast_K}$ is nonnegative, non-increasing and continuous on $(0,\infty)$ and $f^{\ast_K} \leq f^{\ast \ast_K}$.
 \item[$(2')$]
 $$
 (f+g)^{\ast \ast_K}\leq f^{\ast \ast_K}+g^{\ast \ast_K}
 $$
 \item[$(3')$]
 If $|f_n|\uparrow |f|$ $\lambda$-a.e., then $f^{\ast \ast_K}_n\uparrow f^{\ast \ast_K}$.
 \end{itemize}
 For $1\leq p< \infty$ and $1\leq q\leq \infty$, the Lorentz space $L^{p,q} \left( K,\lambda \right)$ is defined as
 $$
 L^{p,q} \left( K,\lambda \right)=\{f: f \,\text{is} \,\lambda \, \text{-measurable on}\, K, \|f\|_{K,p,q}<\infty \}
 $$
 where $\|f\|_{K,p,q}$ is defined by
 $$
 \|f\|_{K,p,q}=\begin{cases}
 \left(\int \limits _0^\infty \left(t^{\frac 1p}f^{\ast \ast_K}(t)\right)^q\frac {dt}{t}\right)^{\frac 1q},&1\leq p<\infty, 1\leq q<\infty\\
 \sup \limits _{t>0}t^{\frac 1p}f^{\ast \ast_K}(t), & 1\leq p\leq \infty, q=\infty.
 \end{cases}
  $$
  Note that if for $1<p\leq\infty$  then $L^{p,p} \left( K,\lambda \right)=L^{p} \left( K,\lambda \right)$. Moreover,
 \begin{equation}\label{LpLpp}
 \|f\|_{K,p}\leq \|f\|_{K,p,p}\leq p'\|f\|_{K,p},
 \end{equation}
  where  $p'=\begin{cases}
\dfrac{p}{p-1},&1< p<\infty, \\
 1, & p= \infty.
 \end{cases} $ \\
  For $p>1$, the space $L^{p,\infty} \left( K,\lambda \right)$ is known as the Marcinkiewicz space or as Weak $L^p\left( K,\lambda \right)$.
  Also note that $L^{1,\infty} \left( K,\lambda \right)=L^{1} \left( K,\lambda \right)$. \\
  If $1<p<\infty$ and $1<q<r<\infty$, then
  $$
  L^{p,q} \left( K,\lambda \right) \subset L^{p,r} \left( K,\lambda \right).
  $$
  Moreover
  \begin{equation} \label{emb}
\|f\|_{K,p,r}\leq \left(\dfrac{q}{p}\right)^{\frac 1q-\frac 1p}\|f\|_{K,p,q}
\end{equation}

The Young inequality on Lebesgue spaces for compact commutative hypergroups was given in \cite{V}. The generalized Young inequality on the Lorentz  spaces for Bessel and Dunkl convolution operators were introduced in \cite{GH} and \cite{HM} correspondingly.

In this paper we establish the generalized Young inequality on the Lorentz  spaces for  commutative hypergroups
and give an application of it to the theory of fractional integrals. The boundedness on the Lorentz spaces of the
fractional integrals on the commutative hypergroups is proved. We also prove the
 Hardy-Littlewood-Sobolev theorem for the fractional integrals on the commutative hypergroups.
\section{Lemmas}
\begin{lemma}
Let $f$ and $\varphi$ be  $\lambda$-measurable functions  on the hypergroup $K$ where $\sup \limits_{x\in K}|f(x)|\leq \beta$ and $f$ vanishes outside of a measurable set $E$ with $\lambda (E)=r$. Then, for $t>0$,
\begin{equation} \label{lemoneil1}
(f\ast_K \varphi)^{\ast \ast_K}(t)\leq \beta r\varphi ^{\ast \ast_K}(r)
\end{equation}
and
\begin{equation} \label{lemoneil2}
(f\ast_K \varphi)^{\ast \ast_K}(t)\leq \beta r\varphi ^{\ast \ast_K}(t).
\end{equation}
\end{lemma}
\begin{Proof}
Without loss of generality we can assume that the functions $f$ and $\varphi$ are nonnegative.
Let $h=f\ast_K\varphi$. For $a>0$, define
$$
 \varphi_{a}(x)=\begin{cases}
 \varphi (x),&\text{if}\,\,\varphi (x) \leq a\\
 a , & \text{if}\,\,\varphi (x) > a,
 \end{cases}
 $$
 $$
 \varphi^{a}(x)=\varphi(x)-\varphi_{a}(x)
 $$
 Also define functions $h_1$ and $h_2$ by
 $$
 h=f\ast_K\varphi _a+f\ast_K\varphi ^a=h_1+h_2.
 $$
 Then we have the following three estimates.
 \begin{equation} \label{h2}
 \sup \limits_{x\in K}h_2(x)\leq \sup \limits_{x\in K}f(x)\|\varphi ^a\|_{K,1}\leq \beta \int\limits_0^\infty \lambda _{\varphi ^a} (s)ds =\beta \int\limits_a^\infty \lambda _\varphi (s)ds,
  \end{equation}
  \begin{equation} \label{h1}
 \sup \limits_{x\in K}h_{1}(x)\leq \|f\|_{K,1}\sup \limits_{x\in K}\varphi _a(x)\leq \beta ra,
  \end{equation}
  and
  \begin{equation} \label{h11}
 \sup \limits_{x\in K}h_2(x)\leq \|f\|_{K,1}\|\varphi ^a\|_{K,1}\leq \beta r \int\limits_a^\infty \lambda _\varphi (s)ds
  \end{equation}
  Now set $a=\varphi ^\ast_K (r)$ in \eqref{h2} and \eqref{h1} and obtain
  $$
  h^{\ast \ast_K}(t)=\frac 1t \int \limits _0^t h^{\ast_K}(s)ds\leq \|h\|_{K,\infty}\leq \|h_1\|_{K,\infty}+\|h_2\|_{K,\infty}
  $$
  $$
  \leq \beta r\varphi ^\ast_K(r)+\beta \int\limits_{\varphi ^\ast_K(r)}^\infty \lambda _\varphi (s)ds,
  $$
  and using \eqref{fstar} we have the inequality \eqref{lemoneil1}.

  Let us prove the inequality \eqref{lemoneil2}. For this purpose set $a=\varphi ^\ast_K (t)$ and use \eqref{h1} and \eqref{h11}. Then
  $$
  th^{\ast \ast_K}(t)=\int \limits _0^t h^{\ast_K}(s)ds\leq\int \limits _0^t h_1^{\ast_K}(s)ds+\int \limits _0^t h_2^{\ast_K}(s)ds
  $$
  $$
  \leq t\|h_1\|_{K,\infty}+\int \limits _0^t h_2^{\ast_K}(s)ds=t\|h_1\|_{K,\infty}+t\|h_2\|_{K,1}
  $$
  $$
  \leq t\beta r\varphi^\ast_K (t)+\beta rt \int\limits_{\varphi ^\ast_K(t)}^\infty \lambda _\varphi (s)ds
  $$
  $$
  =\beta rt\left(\varphi ^\ast_K(t)+ \int\limits_{\varphi ^{\ast_K}(t)}^\infty \lambda _\varphi (s)ds\right)=\beta rt \varphi ^{\ast \ast_K}(t).
  $$
\end{Proof}
\begin{lemma}\label{teoroneil}
Let $f$ and $\varphi$ be $\lambda$-measurable functions on hypergroup $K$, then for all $t>0$ the following inequality holds:
\begin{equation}\label{Oneil}
\left(f\ast_K \varphi\right)^{\ast \ast_K}(t)\leq tf^{\ast \ast_K}(t) \varphi^{\ast \ast_K}(t)+\int \limits_t^\infty f^\ast_K(s)\varphi^\ast_K(s)ds
  \end{equation}
  \end{lemma}
  \begin{Proof}
  Without loss of generality we can assume that the functions $f$ and $\varphi$ are nonnegative. Let $h=f\ast_K \varphi$ and fix $t>0$. Select a nondecreasing sequence $\{s_n\}_{-\infty}^{+\infty}$ such that $s_0=f^\ast_K(t)$, $\lim \limits_{n\rightarrow +\infty}s_n=+\infty$, $\lim \limits_{n\rightarrow -\infty}s_n=0$.\\
  Also let
  $$
  f(x)=\sum \limits_{n=-\infty}^{+\infty}f_n(x)
  $$
  where
  $$
  f_n(x)=\begin{cases}
 0,&\text{if}\,\, f(x)\leq s_{n-1} \\
 f(x)-s_{n-1} & \text{if}\,\,s_{n-1}<f (x) \leq s_n\\
 s_n-s_{n-1} & \text{if}\,\,s_{n}<f (x).
 \end{cases}
  $$
  Since the series $\{s_n\}_{-\infty}^{+\infty}$ converges absolutely we have
  $$
  h= \int \limits_KT^x\varphi(y)\left(\sum \limits_{n=-\infty}^{+\infty}f_n(y^\thicksim)\right)d\lambda (y)
   $$
   $$
   =\sum \limits_{n=-\infty}^{+\infty}\int \limits_KT^x\varphi(y)f_n(y^\thicksim)d\lambda (y)=\sum \limits_{n=-\infty}^{+\infty}\left(f_n\ast_K\varphi\right)
  $$
Define functions $h_1$ and $h_2$ by
  $$
  h=\sum \limits_{n=1}^{+\infty}\left(f_n\ast_K\varphi\right)+\sum \limits_{n=-\infty}^{0}\left(f_n\ast_K\varphi\right)=h_1+h_2.
  $$
  Estimate $h_1^{\ast \ast_K}(t)$. For this purpose use the inequality \eqref{lemoneil2} with $E=\{x:f(x)>s_{n-1}\} $ and $\beta=s_n-s_{n-1}$. We have
  $$
  h_1^{\ast \ast_K}(t)\leq \sum \limits_{n=1}^{+\infty}\left(\left(f_n\ast_K\varphi\right)^{\ast \ast_K}\right)
  $$
  $$
  \leq \sum \limits_{n=1}^{+\infty}(s_n-s_{n-1})\lambda _f(s_{n-1})\varphi ^{\ast \ast_K}(t)
  $$
  $$
  =\varphi ^{\ast \ast_K}(t)\sum \limits_{n=1}^{+\infty}\lambda _f(s_{n-1})(s_n-s_{n-1}).
  $$
  Hence
  \begin{equation} \label{Oneil1}
h_1^{\ast \ast_K}(t) \leq \varphi ^{\ast \ast_K}(t)\int \limits_{f^\ast_K(t)}^\infty \lambda_f(s)ds.
  \end{equation}
  To estimate $h_2^{\ast \ast_K}(t)$ we use the inequality \eqref{lemoneil1}
  $$
  h_2^{\ast \ast_K}(t) \leq \sum \limits_{n=-\infty}^{0}\left(\left(f_n\ast_K\varphi\right)^{\ast \ast_K}\right)
  $$
  $$
  \leq \sum \limits_{n=1}^{+\infty}(s_n-s_{n-1})\lambda _f(s_{n-1})\varphi ^{\ast \ast_K}(\lambda _f(s_{n-1}))
  $$
  $$
  =\sum \limits_{n=1}^{+\infty}\lambda _f(s_{n-1})\varphi ^{\ast \ast_K}(\lambda _f(s_{n-1}))(s_n-s_{n-1}).
  $$
   This implies that
  \begin{equation} \label{Oneil2}
h_2^{\ast \ast_K}(t) \leq \int \limits_0^{f^\ast_K(t)} \lambda_f(s)\varphi ^{\ast \ast_K}(\lambda _f(s))ds.
  \end{equation}
  We will estimate the integral on the right-hand side of \eqref{Oneil2} by making the substitution $s=f^\ast_K (\xi)$ and then integrating by parts. In order to justify the change of variable in the integral, consider a Riemann sum
 $$
  \sum \limits_{n=1}^{+\infty}\lambda _f(s_{n-1})\varphi ^{\ast \ast_K}(\lambda _f(s_{n-1}))(s_n-s_{n-1}).
  $$
  that provides a close approximation to
  $$
  \int \limits_0^{f^\ast_K(t)} \lambda_f(s)\varphi ^{\ast \ast_K}(\lambda _f(s))ds.
  $$
  By adding more points to the Riemann sum if necessary, we may assume that the left-hand end point of each interval on which $\lambda _f$ is constant is included among the $s_n$. then the Riemann sum is not chanched if each $s_n$ that is contained in the interior of an interval on which $\lambda_f$ is constant, is deleted. It is now an easy matter to verify that for each of the remaining $s_n$ there is precisely one element, $\xi_n$, such that $s_n=f^{\ast_K} (\xi_n)$ and that $\lambda \left(f^{\ast_K }(\xi_n)\right)=\xi_n$. Therefore
  $$
  \sum \limits_{n=1}^{+\infty}\lambda _f(s_{n-1})\varphi ^{\ast \ast_K}(\lambda _f(s_{n-1}))(s_n-s_{n-1}).
  $$
  $$
  =\sum \limits_{n=1}^{+\infty}\xi_{n-1}\varphi ^{\ast \ast_K}(\xi_{n-1})(f^{\ast_K} (\xi_n)-f^{\ast_K} (\xi_{n-1}))
  $$
  which, by adding more points if necessary, provides a close approximation to
  $$
  -\int \limits_t^\infty
  \xi\varphi ^{\ast \ast_K}(\xi)df^{\ast_K} (\xi).
  $$
  If we recall \eqref{Oneil2} we get
  \begin{equation} \label{Oneil3}
h_2^{\ast \ast_K}(t) \leq \int \limits_0^{f^{\ast_K}(t)} \lambda_f(s)\varphi ^{\ast \ast_K}(\lambda _f(s))ds=-\int \limits_t^\infty \xi\varphi ^{\ast \ast_K}(\xi)df^{\ast_K} (\xi).
  \end{equation}
  Now let $\delta$ be an arbitrarily large number and choose $\xi_j$ such that $t=\xi_1\leq \xi_2\leq \ldots \leq \xi_{j+1}=\delta$. Then
  $$
  \delta\varphi ^{\ast \ast_K}(\delta)f^{\ast_K}(\delta) -t\varphi ^{\ast \ast_K}(t)f^{\ast_K}(t)
  $$
  $$
  =\sum \limits _{n=1}^j \xi_{n+1}\varphi ^{\ast \ast_K}(\xi_{n+1})\left(f^{\ast_K}(\xi_{n+1})-f^{\ast_K}(\xi_{n})\right)
$$
$$
  +\sum \limits _{n=1}^j f^{\ast_K}(\xi_{n})\left(\varphi ^{\ast \ast_K}(\xi_{n+1})\xi_{n+1}-\varphi ^{\ast \ast_K}(\xi_{n})\xi_{n}\right)
$$
$$
  =\sum \limits _{n=1}^j \xi_{n+1}\varphi ^{\ast \ast_K}(\xi_{n+1})\left(f^{\ast_K}(\xi_{n+1})-f^{\ast_K}(\xi_{n})\right)
$$
$$
  +\sum \limits _{n=1}^j f^\ast_K(\xi_{n})\int \limits_{\xi_{n}}^{\xi_{n+1}} \varphi ^{\ast_K}(\tau)d\tau
$$
$$
  \leq \sum \limits _{n=1}^j \xi_{n+1}\varphi ^{\ast \ast_K}(\xi_{n+1})\left(f^{\ast_K}(\xi_{n+1})-f^{\ast_K}(\xi_{n})\right)
$$
$$
+\sum \limits _{n=1}^j f^{\ast_K}(\xi_{n})\varphi ^{\ast_K}(\xi_{n})\left(\xi_{n+1}-\xi_{n}\right).
$$
This means that
\begin{equation} \label{Oneil4}
\delta\varphi ^{\ast \ast_K}(\delta)f^{\ast_K}(\delta) -t\varphi ^{\ast \ast_K}(t)f^{\ast_K}(t)\leq \int \limits_t^\delta \xi\varphi ^{\ast \ast_K}(\xi)df^{\ast_K} (\xi)+\int \limits_t^\delta f^{\ast_K}(\xi)\varphi ^{\ast_K}\xi)d\xi  .
  \end{equation}
  Now we estimate the expression $\delta\varphi ^{\ast \ast_K}(\delta)f^{\ast_K}(\delta) -t\varphi ^{\ast \ast_K}(t)f^{\ast_K}(t)$ below.
  $$
  \delta\varphi ^{\ast \ast_K}(\delta)f^{\ast_K}(\delta) -t\varphi ^{\ast \ast_K}(t)f^{\ast_K}(t)
  $$
  $$
  =\sum \limits _{n=1}^j \xi_{n}\varphi ^{\ast \ast_K}(\xi_{n})\left(f^{\ast_K}(\xi_{n+1})-f^{\ast_K}(\xi_{n})\right)
$$
$$
  +\sum \limits _{n=1}^j f^{\ast_K}(\xi_{n+1})\left(\varphi ^{\ast \ast_K}(\xi_{n+1})\xi_{n+1}-\varphi ^{\ast \ast_K}(\xi_{n})\xi_{n}\right)
$$
$$
  =\sum \limits _{n=1}^j \xi_{n}\varphi ^{\ast \ast_K}(\xi_{n})\left(f^{\ast_K}(\xi_{n+1})-f^{\ast_K}(\xi_{n})\right)
$$
$$
  +\sum \limits _{n=1}^j f^{\ast_K}(\xi_{n+1})\int \limits_{\xi_{n}}^{\xi_{n+1}} \varphi ^{\ast_K}(\tau)d\tau
$$
$$
  \geq \sum \limits _{n=1}^j \xi_{n}\varphi ^{\ast \ast_K}(\xi_{n})\left(f^{\ast_K(}\xi_{n+1})-f^{\ast_K}(\xi_{n})\right)
$$
$$
+\sum \limits _{n=1}^j f^\ast_K(\xi_{n+1})\varphi ^{\ast_K}(\xi_{n+1})\left(\xi_{n+1}-\xi_{n}\right).
$$
In other words
\begin{equation} \label{Oneil5}
\delta\varphi ^{\ast \ast_K}(\delta)f^{\ast_K}(\delta) -t\varphi ^{\ast \ast_K}(t)f^{\ast_K}(t)\geq \int \limits_t^\delta \xi\varphi ^{\ast \ast_K}(\xi)d f^{\ast_K }(\xi)+\int \limits_t^\delta f^{\ast_K}(\xi)\varphi ^{\ast_K}(\xi)d\xi  .
  \end{equation}
  From \eqref{Oneil4} and \eqref{Oneil5} we obtain
  $$
   -\int \limits_t^\delta \xi\varphi ^{\ast \ast_K}(\xi)df^{\ast_K} (\xi)=t\varphi ^{\ast \ast_K}(t)f^{\ast_K}(t)-\delta\varphi ^{\ast \ast_K}(\delta)f^{\ast_K}(\delta) +\int \limits_t^\delta f^{\ast_K}(\xi)\varphi ^{\ast_K}(\xi)d\xi  .
  $$
  $$
  \leq t\varphi ^{\ast \ast_K}(t)f^{\ast_K}(t)+\int \limits_t^\delta f^{\ast_K}(\xi)\varphi ^{\ast_K}(\xi)d\xi.
  $$
  Thus
  $$
   -\int \limits_t^\infty \xi\varphi ^{\ast \ast_K}(\xi)df^{\ast_K }(\xi)
  \leq t\varphi ^{\ast \ast_K}(t)f^{\ast_K}(t)+\int \limits_t^\infty f^{\ast_K}(\xi)\varphi ^{\ast_K}(\xi)d\xi.
  $$
  By using this inequality and \eqref{Oneil3} we have
  \begin{equation} \label{Oneil6}
h_2^{\ast \ast_K}(t) \leq \int \limits_0^{f^{{\ast_K}}(t)} \lambda_f(s)\varphi ^{\ast \ast_K}(\lambda _f(s))ds\leq t\varphi ^{\ast \ast_K}(t)f^{\ast_K}(t)+\int \limits_t^\infty f^{\ast_K}(\xi)\varphi ^{\ast_K}(\xi)d\xi.
  \end{equation}
  Finally, from \eqref{Oneil1}, \eqref{Oneil6} and \eqref{fstar} we get
  $$
 h^{\ast \ast_K}(t)\leq h_1^{\ast \ast_K}(t)+ h_2^{\ast \ast_K}(t)
 $$
 $$
 \leq \varphi ^{\ast \ast_K}(t)\int \limits_{f^{\ast_K}(t)}^\infty \lambda_f(s)ds+t\varphi ^{\ast \ast_K}(t)f^{\ast_K}(t)+\int \limits_t^\infty f^{\ast_K}(\xi)\varphi ^{\ast_K}(\xi)d\xi
  $$
  $$
 = f^{\ast_K}(t)\varphi ^{\ast \ast_K}(t)+\int \limits_t^\infty f^{\ast_K}(\xi)\varphi ^{\ast_K}(\xi)d\xi
  $$
  $$
 = tf^{\ast \ast_K}(t)\varphi ^{\ast \ast_K}(t)+\int \limits_t^\infty f^{\ast_K}(\xi)\varphi ^{\ast_K}(\xi)d\xi.
  $$
\end{Proof}
\begin{lemma}
Let $f$ and $\varphi$ be $\lambda$-measurable functions on hypergroup $K$, then for all $t>0$ the following inequality holds:
\begin{equation} \label{2Oneil}
\left(f\ast_K \varphi\right)^{\ast \ast_K}(t)\leq \int \limits_t^\infty f^{\ast \ast_K}(s)\varphi^{\ast \ast_K}(s)ds
  \end{equation}
  \end{lemma}
  \begin{Proof}
  Assume that the integral on the right of \eqref{2Oneil} is finite. Then it is easy to see
  \begin{equation} \label{2Oneil1}
 sf^{\ast \ast_K}(s)\varphi^{\ast \ast_K}(s)\rightarrow 0, \,\text{as}\, s\rightarrow \infty .
  \end{equation}
  Let $h=f\ast_K\varphi$.\\
  By Lemma \ref{teoroneil} we have
  $$
  h^{\ast \ast_K}(t)\leq tf^{\ast \ast_K}(t) \varphi^{\ast \ast_K}(t)+\int \limits_t^\infty f^{\ast_K}(s)\varphi^{\ast_K}(s)ds
  $$
  \begin{equation} \label{2Oneil2}
  \leq tf^{\ast \ast_K}(t) \varphi^{\ast \ast_K}(t)+\int \limits_t^\infty f^{\ast \ast_K}(s)\varphi^\ast_K(s)ds.
  \end{equation}
  Since $f^{\ast \ast_K}$ and $g^{\ast \ast_K}$ are non-increasing,
  \begin{equation} \label{2Oneil3}
  \frac{df^{\ast \ast_K}(s)}{ds}=-\frac{1}{s^2}\int \limits _0^sf^{\ast_K}(\tau)d\tau+\frac{1}{s}f^{\ast_K}(s)=\frac{1}{s}\left(f^{\ast_K}(s)-f^{\ast \ast_K}(s)\right),
  \end{equation}
  \begin{equation} \label{2Oneil4}
  \frac{d(s\varphi^{\ast \ast_K}(s))}{ds}=\varphi^{\ast \ast_K}(s)+s\left(\frac{1}{s}\left(\varphi^{\ast_K}(s)
  -\varphi^{\ast \ast_K}(s)\right)\right)=\varphi^{\ast_K}(s)
  \end{equation}
  for $m$-almost all $s$. Since $f^{\ast \ast_K}$ and $g^{\ast \ast_K}$ are absolutely continuous, we may use the integration by parts for $\int \limits_t^\infty f^{\ast \ast_K}(s)d\left(s\varphi^{\ast \ast_K}(s)\right)$. Using \eqref{2Oneil3}, \eqref{2Oneil4} and \eqref{2Oneil1} we obtain
    $$
  \int \limits_t^\infty f^{\ast \ast_K}(s)\varphi^{\ast_K}(s)ds
   =\int \limits_t^\infty f^{\ast \ast_K}(s)d\left(s\varphi^{\ast \ast_K}(s)\right)
  $$
  $$
  =\left.f^{\ast \ast_K}(s)s\varphi^{\ast \ast_K}(s)\right|_t^\infty-\int \limits_t^\infty s\varphi^{\ast \ast_K}(s)df^{\ast \ast_K}(s)
   $$
   $$
   =-tf^{\ast \ast_K}(t)\varphi^{\ast \ast_K}(t)+\int \limits_t^\infty \varphi^{\ast \ast_K}(s)(f^{\ast \ast_K}(s)-f^{\ast_K}(s))ds
   $$
   \begin{equation} \label{2Oneil5}
   \leq -tf^{\ast \ast_K}(t)\varphi^{\ast \ast_K}(t)+\int \limits_t^\infty \varphi^{\ast \ast_K}(s)f^{\ast \ast_K}(s)ds
   \end{equation}
   By \eqref{2Oneil2} and \eqref{2Oneil5} we have
   $$
  h^{\ast \ast_K}(t)\leq
  \int \limits_t^\infty f^{\ast \ast_K}(s)\varphi^{\ast \ast_K}(s)ds.
  $$
  \end{Proof}
  The next lemma is a classical estimate, known as Hardy's inequality.
  \begin{lemma}$\left(\cite{BSh}, \cite{Z}\right)$
  If  $1\leq p<\infty$, $q>0$ and $f$ be a nonnegative $m$-measurable function on $(0,\infty)$, then
  \begin{equation} \label{hardy}
  \int \limits_0^\infty \left(\frac 1s \int \limits_0^s f(\tau)d\tau\right)^ps^{p-q-1}ds\leq \left(\frac{p}{q}\right)^q\int \limits_0^\infty f(t)^pt^{p-q-1}dt.
  \end{equation}
  \end{lemma}
  \section{Generalization of Young's inequality}
  \begin{theorem}\label{young}
  If  $f\in L^{p_1,q_1} \left( K,\lambda \right)$, $\varphi \in L^{p_2,q_2} \left( K,\lambda \right)$ and $\dfrac{1}{p_1}+\dfrac{1}{p_2}>1$, then $(f\ast_K\varphi)\in L^{p_0,q_0} \left( K,\lambda \right)$ where $\dfrac{1}{p_1}+\dfrac{1}{p_2}-1=\dfrac{1}{p_0}$ and $q_0\geq 1$ is any number such that
  $\dfrac{1}{q_1}+\dfrac{1}{q_2}\geq \dfrac{1}{q_0}$.\\
  Moreover,
  \begin{equation}\label{lorentz}
\|(f\ast_K\varphi)\|_{K,p_0,q_0}\leq 3p_0\|f\|_{K,p_1,q_1}\|\varphi\|_{K,p_2,q_2}.
\end{equation}
    \end{theorem}
  \begin{Proof}
  Let $h=f\ast_K\varphi$.\\
  Suppose that $q_1$, $q_2$, $q_0$ are all different from $\infty$. Then, by  \eqref{2Oneil}, we have
  $$
  (\|h\|_{K,p_0,q_0})^{q_0}=\int \limits _0^\infty \left(s^{\frac{1}{p_0}}h^{\ast \ast_K}(s)\right)^q \frac{ds}{s}
  $$
  $$
  \leq\int \limits _0^\infty \left(s^{\frac{1}{p_0}}\int \limits _s^\infty f^{\ast \ast_K}(\tau)\varphi^{\ast \ast_K}(\tau)d\tau\right)^q \frac{ds}{s}
  $$
  $$
   =\int \limits _0^\infty \left(\frac{1}{t^{\frac{1}{p_0}}}\int \limits _0^t f^{\ast \ast_K}\left(\frac{1}{\eta}\right)\varphi^{\ast \ast_K}\left(\frac{1}{\eta}\right)\frac{d\eta}{\eta ^2}\right)^q \frac{dt}{t}.
  $$
  The last equality was obtained by the change of variables $s=\dfrac{1}{t}$ and $\tau=\dfrac{1}{\eta}$. Using \eqref{hardy} we get
  $$
   \int \limits _0^\infty \left(\frac{1}{t^{\frac{1}{p_0}}}\int \limits _0^t f^{\ast \ast_K}\left(\frac{1}{\eta}\right)\varphi^{\ast \ast_K}\left(\frac{1}{\eta}\right)\frac{d\eta}{\eta ^2}\right)^q \frac{dt}{t}
  $$
  $$
   \leq p_{0}^{q_0}\int \limits _0^\infty \left(t^{1-\frac {1}{p_0}} \frac {f^{\ast \ast_K}\left(\frac{1}{t}\right)\varphi^{\ast \ast_K}\left(\frac{1}{t}\right)}{t^2}\right)^{q_0}\frac{dt}{t }
  $$
  $$
    =p_{0}^{q_0}\int \limits _0^\infty \left(s^{1+\frac {1}{p_0}} f^{\ast \ast_K}\left(s\right)\varphi^{\ast \ast_K}\left(s\right)\right)^{q_0}\frac{ds}{s}
  $$
  The last equality was obtained by the change of the variable $t=\dfrac{1}{s}$. Since $\dfrac{q_0}{q_1}+\dfrac{q_0}{q_2}\geq 1$, one can find positive numbers $n_1$ and $n_2$ such that
  $$
  \frac{1}{n_1}+\frac{1}{n_2}=1 \,\, \text{and}\,\, \frac{1}{n_1}\leq \frac{q_0}{q_1},\,\, \frac{1}{n_2}\leq \frac{q_0}{q_2}.
  $$
  By H$\ddot{\text{o}}$lder's inequality we obtain
  $$
   (\|h\|_{K,p_0,q_0})^{q_0} \leq p_{0}^{q_0}\int \limits _0^\infty \frac {\left(s^{\frac {1}{p_1}}f^{\ast \ast_K}\left(s\right)\right)^{q_0}}{s^{\frac{1}{n_2}}}\frac {\left(s^{\frac {1}{p_2}}\varphi^{\ast \ast_K}\left(s\right)\right)^{q_0}}{s^{\frac{1}{n_1}}}ds
  $$
  $$
    \leq p_{0}^{q_0}\left[\int \limits _0^\infty \left(s^{\frac {1}{p_1}}f^{\ast \ast_K}\left(s\right)\right)^{q_0n_1}\frac {ds}{s}\right]^{\frac {1}{n_1}}\left[\int \limits _0^\infty \left(s^{\frac {1}{p_2}}\varphi^{\ast \ast_K}\left(s\right)\right)^{q_0n_2}\frac {ds}{s}\right]^{\frac {1}{n_2}}
  $$
  $$
  =p_{0}^{q_0}\left(\|f\|_{K,p_1,q_0n_1}\right)^{q_0}\left(\|\varphi\|_{K,p_2,q_0n_2}\right)^{q_0}.
  $$
 Finally, by \eqref{emb} we have
 $$
\|h\|_{K,p_0,q_0} \leq p_{0}\|f\|_{K,p_1,q_0n_1} \|\varphi\|_{K,p_2,q_0n_2}\leq p_0e^{\frac 1e}e^{\frac 1e}\|f\|_{K,p_1,q_1} \|\varphi\|_{K,p_2,q_2}\leq 3p_0\|f\|_{K,p_1,q_1} \|\varphi\|_{K,p_2,q_2}.
$$
Similar reasoning leads to the desired result in case one or more of $q_1$, $q_2$, $q_0$ are $\infty$.
    \end{Proof}
    \section{Applications to the theory of  fractional integrals  }
Consider the following particular case of  Theorem \ref{young}. If we take $p_1=\dfrac{N}{N-\alpha}$, with $o<\alpha<N$, $q_1=\infty$ in Theorem \ref{young}, then the condition $\dfrac {1}{p_1}+\dfrac {1}{p_2}>1$ is equivalent to $\alpha<\dfrac{N}{p_2}$, and the condition $\dfrac {1}{p_1}+\dfrac {1}{p_2}-1=\dfrac{1}{p_0}$ is equivalent to $\dfrac{1}{p_0}=\dfrac{1}{p_2}-\dfrac{\alpha}{N}$. Thus we have the following result.
\begin{theorem}\label{fraclorentz}
Let $(K,\ast_K)$ be a commutative hypergroup, with  Haar measure $\lambda$.
  If  $f\in L^{\frac{N}{N-\alpha},\infty} \left( K,\lambda \right)$, $\varphi \in L^{p,q} \left( K,\lambda \right)$, where $0<\alpha<\dfrac{N}{p}$, $1\leq q\leq \infty$  then $(f\ast_K\varphi)\in L^{r,q} \left( K,\lambda \right)$
  and
  \begin{equation}\label{kalfaineq}
\|(f\ast_K\varphi)\|_{K,r,q}\leq 3r\|f\|_{K,\frac{N}{N-\alpha},\infty}\|\varphi\|_{K,p,q}.
\overline{}\end{equation}
  where $\dfrac{1}{r}=\dfrac{1}{p}-\dfrac{\alpha}{N}$.
    \end{theorem}
    Let $K$ be a set. A function $\rho
:K\times K\rightarrow \left[ 0,\infty \right) $ is called
quasi-metric if:

\begin{enumerate}
\item
$\rho \left( x,y\right) =0\,\Leftrightarrow \;x=y ;$
\item
$\rho \left( x,y\right) =\rho \left( y,x\right) ;$
\item
there is a constant $c\geq 1$ such that for every $x,y,z\in X$
$$
\rho \left( x,y\right) \leq c\left( \rho \left( x,z\right) +\rho
\left( z,y\right) \right) .
$$
\end{enumerate}
    Define  the fractional integral (or Riesz potential)
$$
I_\alpha f(x)=\int \limits_KT^x\rho(e,y)^{\alpha -N}f(y^\thicksim)d\lambda(y), \,\,0<\alpha<N
$$
on commutative hypergroup $(K,\ast_K)$ equipped with the pseudo-metric $\rho$.\\
Also define a ball $B(e,r)=\{y\in K: \, \rho(e,y)<r\}$ with a center $e$ and a radius $r$.
   \begin{theorem}Let $(K,\ast_K)$ be a commutative hypergroup, with quasi-metric $\rho$ and  Haar measure $\lambda$ satisfying $\lambda B(e,r)= Ar^N$, where $A$ is a positive constant.
      Assume that   $1\leq q\leq \infty$, $1\leq p< \infty$, $0<\alpha<\dfrac{N}{p}$.  If $f \in L^{p,q} \left( K,\lambda \right)$,   then $I_\alpha f\in L^{r,q} \left( K,\lambda \right)$
  and
  \begin{equation}\label{kalfaineq}
\|I_\alpha f\|_{K,r,q}\leq C\|f\|_{K,p,q},
\end{equation}
  where $\dfrac{1}{r}=\dfrac{1}{p}-\dfrac{\alpha}{N}$ and  $C=\dfrac {3rN}\alpha A^{\frac{N-\alpha}{N}}$.
    \end{theorem}
    \begin{Proof}
    Let us show that $\rho(e,\cdot)^{\alpha -N} \in L^{\frac{N}{N-\alpha},\infty}\left( K,\lambda \right)$. For the distribution of $\rho(e,\cdot)^{\alpha -N}$ we can write
    $$
    \lambda_{\rho(e,\cdot)^{\alpha -N}}(t)= \lambda\{x: x\in K, \rho(e,x)^{\alpha -N}>t\}
    $$
    $$
   = \lambda\{x: x\in K, \rho(e,x)<t^{\frac{1}{\alpha-N}}\}=At^{\frac{N}{\alpha-N}}.
    $$
    Since $\rho(e,\cdot)^{\alpha -N}$ is continuous and strictly decreasing we have $\left(\rho(e,\cdot)^{\alpha -N}\right)^{\ast_K}$ is the inverse of the distribution function. That is  $\left(\rho(e,\cdot)^{\alpha -N}\right)^{\ast_K}(t)=\left(\dfrac{t}{A}\right)^{\frac{\alpha-N}{N}}$. Then
    $$
    \left(\rho(e,\cdot)^{\alpha -N}\right)^{\ast \ast_K}(t)=\frac 1t \int \limits_{0}^{t}\left(\dfrac{s}{A}\right)^{\frac{\alpha-N}{N}}ds=\dfrac N\alpha \left(\dfrac tA\right)^{\frac{\alpha-N}{N}}.
    $$
    Therefore $\rho(e,\cdot)^{\alpha -N} \in L^{\frac{N}{N-\alpha},\infty}\left( K,\lambda \right)$ and
    \begin{equation} \label{ialfa}
    \|\rho(e,\cdot)^{\alpha -N} \|_{K,\frac{N}{N-\alpha},\infty}=\dfrac N\alpha A^{\frac{N-\alpha}{N}}.
    \end{equation}
    Thus, from Theorem \ref{fraclorentz} and \label{ialfa} we have the required result.
    \end{Proof}
 \begin{theorem} \label{fraclebeq}
 Let $(K,\ast_K)$ be a commutative hypergroup, with  Haar measure $\lambda$.
 If  $f\in L^{\frac{N}{N-\alpha},\infty} \left( K,\lambda \right)$, $\varphi \in L^{p} \left( K,\lambda \right)$, where $1<p< \infty$, $0<\alpha<\dfrac{N}{p}$,  then $(f\ast_K\varphi)\in L^{r} \left( K,\lambda \right)$
  and
  \begin{equation}\label{kalfaineq2}
\|(f\ast_K\varphi)\|_{K,r}\leq 3r\dfrac{p}{p-1}\left(\dfrac{p}{r}\right)^{\frac 1p-\frac 1r}\|f\|_{K,\frac{N}{N-\alpha},\infty}\|\varphi\|_{K,p}.
\overline{}\end{equation}
  where $\dfrac{1}{r}=\dfrac{1}{p}-\dfrac{\alpha}{N}$.
  \end{theorem}
  \begin{Proof} From \eqref{LpLpp}, \eqref{emb} and \eqref{lorentz} we have
  $$
  \|(f\ast_K\varphi)\|_{K,r}\leq \|(f\ast_K\varphi)\|_{K,r,r}\leq \left(\dfrac{p}{r}\right)^{\frac 1p-\frac 1r}\|(f\ast_K\varphi)\|_{K,r,p}
  $$
  $$
  \leq 3r\left(\dfrac{p}{r}\right)^{\frac 1p-\frac 1r}\|f\|_{K,\frac{N}{N-\alpha},\infty}\|\varphi\|_{K,p,p}\leq 3r\dfrac{p}{p-1}\left(\dfrac{p}{r}\right)^{\frac 1p-\frac 1r}\|f\|_{K,\frac{N}{N-\alpha},\infty}\|\varphi\|_{K,p}
  $$
  \end{Proof}
  The following result give us the Hardy-Littlewood-Sobolev theorem for the fractional integrals on the commutative hypergroups.
  \begin{theorem}Let $(K,\ast_K)$ be a commutative hypergroup, with quasi-metric $\rho$ and  Haar measure $\lambda$ satisfying $\lambda B(e,r)= Ar^N$, where $A$ is a positive constant.
      Assume that   $1 <p< \infty$, $0<\alpha<\dfrac{N}{p}$.  If $f \in L^{p} \left( K,\lambda \right)$,   then $I_\alpha f\in L^{r} \left( K,\lambda \right)$
  and
  \begin{equation}\label{kalfaineq}
\|I_\alpha f\|_{K,r}\leq C\|f\|_{K,p},
\end{equation}
  where $\dfrac{1}{r}=\dfrac{1}{p}-\dfrac{\alpha}{N}$ and  $C=\dfrac{3pr}{p-1}\left(\dfrac{p}{r}\right)^{\frac 1p-\frac 1r}\dfrac {N}\alpha A^{\frac{N-\alpha}{N}}$.
    \end{theorem}
    \begin{Proof} This follows immediately from Theorem \ref{fraclebeq} and \eqref{ialfa}.
    \end{Proof}

    \textbf{Acknowledgement.} This work was supported by the Science Development Foundation
under the President of the Republic of Azerbaijan Grant EIF-2012-2(6)-39/10/1. The author
would like to express his thanks to Academician Akif Gadjiev for valuable remarks.

Mubariz G. Hajibayov

National Aviation Academy.
Bine gesebesi, 25-ci km, AZ1104, Baku, Azerbaijan

and

Institute of Mathematics and Mechanics of NAS of Azerbaijan.%

9, B. Vahavzade str., AZ1141, Baku, Azerbaijan.

\end{document}